%% file: wdgrii.tex
\documentclass[11pt,leqno]{article}

\input ./header.tex

\begin{document}
\title{A Class of Representations of Hecke Algebras II}
\author{Dean Alvis}

\date{}

\maketitle

\begin{abstract}
Let $W$ be a Coxeter group whose proper parabolic 
subgroups are finite.
According to 
Theorem~1.12 of \cite{digraphpaper},
 if the module of
a finite $W$-digraph $\Gamma$ is isomorphic to the module
of a $W$-graph over $\rationals$, then $\Gamma$ is acyclic.  
We extend this result to 
 Coxeter groups with finite dihedral 
 parabolic subgroups and 
$W$-graphs over arbitrary fields $F \le \complexes$. 
Also, an example is provided showing the converse of this 
theorem is false. That is, there is an example
of a finite, acyclic
$W$-digraph whose module does not afford a $W$-graph.
\end{abstract}


\section{An extension of Theorem~1.12 of  \cite{digraphpaper}}

Let $(W,S)$ be a Coxeter system with presentation
\[
W =
\left<
s \in S
\, \mid \,
(rs)^{n(r,s)}=e \text{ for $r,s\in S$ whenever $n(r,s)<\infinity$}
\right> ,
\]
where $n(s,s) = 1$ and $1<n(r,s)=n(s,r)\le\infinity$ for $r,s\in S$, $r \ne s$.
Let $\ell$ be the length function of $(W,S)$.
Let $u$ be an indeterminate over $\complexes$, 
and let $H$ be the Hecke algebra
of $(W,S)$ over $\rationals(u)$.  

See
\cite{digraphpaper} for the definition of {\it $W$-digraph}.
For $\Gamma$ a $W$-digraph and $\beta \in \vertices(\Gamma)$, 
define $\In(\beta)$ to be the set of all $s \in S$ such that
$\Gamma$ has an edge of the form
\solidedge{$\alpha$}{$\beta$}{$s$}
or
\dashededge{$\alpha$}{$\beta$}{$s$} 
for some $\alpha \in \vertices(\Gamma)$. 
Observe $\beta$ is a source (sink) in $\Gamma$ if and only if 
$\In(\beta) = \emptyset$ ($\In(\beta) = S$, respectively). 
For $J \subseteq S$, put
\[
N_{\Gamma}(J) =
\abs{ \setof{\beta \in \vertices(\Gamma)} {\In(\beta) = J} }.
\]
(This definition will only be applied when $\Gamma$ is finite, i.e. 
when $\vertices(\Gamma)$ and $\edges(\Gamma)$ are finite.)

Let $\Psi$ be a $W$-graph over the subfield $F$ of $\complexes$, 
in the sense of \cite{gyojawgraph}, Definition~2.1, with $u$ playing the role of
 $q^{1/2}$ in \cite{gyojawgraph}.  Thus $I$ is the
vertex-labeling function $x \mapsto I_x \subseteq S$, $x \in X$, 
and $\mu$ is the 
edge-labeling function $\mu : X \times X \setminus \Delta  \rightarrow F$.
For  $J \subseteq S$, put
\[
N_{\Psi}(J) =
\abs{ \setof{x \in \vertices(\Psi)} { I_x = J} }.
\]
(This definition will only be applied when $\Psi$ is finite, i.e. 
when $\vertices(\Psi)$ is finite.)


The goal of this section is to prove the following.

\begin{theorem}
If $n(s,t) < \infinity$ for $s, t \in S$, 
$\Gamma$ is a finite $W$-digraph, $\Psi$ is a $W$-graph
over a subfield $F$ of $\complexes$, 
and $M(\Gamma)^F = F(u) \otimes_{\rationals(u)} M(\Gamma)$ is 
isomorphic to $M(\Psi)$ as $H^F$-modules, then
the following hold:
\begin{enumerate}[{\upshape(i)}]
\item
$N_{\Gamma}(J) = N_{\Psi}(J)$ for all $J \subseteq S$.
\item
$\Gamma$ is acyclic.
\end{enumerate}
\label{theorem:wgraphacyclic}
\end{theorem}


We require the following results.
For $F$ a subfield of $\complexes$, $\lambda$ a linear character
of $H^F$, and $V$ an $H^F$-module, define 
\[
V_\lambda =
\setof{v \in V}{h v = \lambda(h) v \text{ for all } h \in H^F}
\]
and 
\[
\charprod{V}{\lambda}{H^F}
=
\dim_{F(u)} V_{\lambda}.
\]
Also, let $\ind^F$ and $\sgn^F$ be the linear characters of $H^F$ determined
by $\ind^F(T_w) = u_w = u^{2 \ell(w)}$ and
$\sgn^F(T_w) = \varepsilon_w = (-1)^{\ell(w)}$ for $w \in W$.


\begin{lemma}
\label{lemma:eigenvalues}
Let $F$ a subfield of $\complexes$. Let $\Gamma$ be
a $W$-digraph, and put $V = M(\Gamma)^F$.
Suppose $s \in S$ and
 $v = \sum_{\gamma\in X} \lambda_\gamma \gamma \in V$.
Then the following hold:
\begin{enumerate}[{\upshape(i)}]
\item
$T_s v = u^2 v$ if and only if
$\lambda_\beta = \lambda_\alpha$ whenever 
\solidedge{$\alpha$}{$\beta$}{$s$}
or
\dashededge{$\alpha$}{$\beta$}{$s$}
is an edge of $\Gamma$.
\item
$T_s v = - v$ if and only if
\[
\lambda_\beta
=
\begin{cases}
- u^{-2} \lambda_\alpha & 
\text{whenever $\solidedge{$\alpha$}{$\beta$}{$s$} \in \edges(\Gamma)$,} \\
- (u+1)(u^2-u)^{-1} \lambda_\alpha &
\text{whenever $\dashededge{$\alpha$}{$\beta$}{$s$} \in \edges(\Gamma)$.}
\end{cases}
\]
\end{enumerate}
\end{lemma}

\begin{proof}
The argument given for 
Lemma~2.4 in \cite{digraphpaper} applies, 
with $F$, $H^F$, and $V=M(\Gamma)^F$ replacing
$\rationals$, $H$, and $M=M(\Gamma)$, respectively.
\end{proof}


\begin{lemma}
\label{lemma:linearcharmults}
If $\Gamma$ is a $W$-digraph, $\vertices(\Gamma)$ is finite, 
and $F$ is a subfield of $\complexes$, 
then the following hold:
\begin{enumerate}[{\upshape(i)}]
\item
The number of connected components of $\Gamma$ is equal to 
$\charprod{M(\Gamma)^F}{\ind^F}{H^F}$.
\item
If $n(s,t) < \infinity$ for all $s, t \in S$, then
the number of acyclic connected components of $\Gamma$ is equal
to $\charprod{M(\Gamma)^F}{\sgn^F}{H^F}$.
\end{enumerate}
\end{lemma}

\begin{proof}
The proof for  Theorem~1.7 in \cite{digraphpaper} 
 applies here if
$\rationals$, $H$, and $M(\Gamma)$ are replaced by
by $F$, $H^F$, $M(\Gamma)^F$, respectively, using 
Lemma~\ref{lemma:eigenvalues} in place of  
\cite{digraphpaper}, Lemma~2.4.
\end{proof}


\begin{lemma}
If $F$ is a subfield of $\complexes$, $\Psi$ is a $W$-graph over $F$,  
$x \mapsto I_x \subseteq S$ is the vertex-labeling function for $\Psi$, 
 and $M(\Psi)_{\ind^F} \ne \set{0}$, 
then there is some $x_0 \in \vertices(\Psi)$ such that $I_{x_0}= \emptyset$.
\label{lemma:indlemma}
\end{lemma}

\begin{proof}
Let $X = \vertices(\Psi)$, and let 
$\mu$ be the edge-labeling function of $\Psi$. 
Suppose
$v = \sum_{x \in X} \gamma_x x \in M(\Psi)_{\ind^F}$ and $v \ne 0$, where
 all but finitely many coefficients $\gamma_x$ are zero.
Replacing $v$ by a scalar multiple if necessary, we can assume 
$\gamma_x \in F[u]$ for all $x\in X$ and
$\gcd \setof{\gamma_x}{x \in X} = 1$.
Choose $x_0 \in X$ such that $\gamma_{x_0} \not\in u F[u]$.
Suppose $I_{x_0} \ne \emptyset$, so that
 $s \in I_{x_0}$ for some $s \in S$.  Since 
\[
u^2 v = T_s v
 = 
- \sum_{ x \in X, s \in I_x } \gamma_x x
+
\sum_{ y \in X, s \not\in I_y }
 \gamma_y \left( u^2 y + u \sum_{z \in X, s \in I_z} \mu(z,x) z \right),
\]
comparing coefficients of $x_0$ shows $\gamma_{x_0} \in u F[u]$, and
so a contradiction is reached.  
Therefore $I_{x_0}= \emptyset$.
\end{proof}


\begin{proof}[Proof of Theorem~\ref{theorem:wgraphacyclic}]
For the remainder of this proof it is
assumed that  $n(s,t) < \infinity$ for $s, t \in S$, 
$\Gamma$ is a finite $W$-digraph, $\Psi$ is a $W$-graph
over the subfield $F$ of $\complexes$, 
and $M(\Gamma)^F$ is isomorphic to $M(\Psi)$ as $H^F$-modules.

Any connected component of $\Gamma$ that
contains a sink is acyclic by 
\cite{digraphpaper}, Theorem~1.5(ii).  On the other hand,
any acyclic connected component $C$ of 
$\Gamma$ contains some sink $\sigma$ because  
$\Gamma$ is finite, 
and $\sigma$ is the unique sink
in $C$ by \cite{digraphpaper}, Theorem~1.5(i).
Thus the number of sinks in $\Gamma$, that is, 
$N_\Gamma(S)$, is equal
to the number of acyclic connected components of $\Gamma$.
Thus by  
Lemma~\ref{lemma:linearcharmults}(ii),
$N_\Gamma(S)$ 
is equal to 
$\charprod{M(\Gamma)^F}{\sgn^F}{H^F}$. 
Also, 
$\charprod{M(\Psi)}{\sgn^F}{H^F} = N_{\Psi}(S)$
because
$M(\Psi)_{\sgn^F}$ has basis 
$\setof{x \in \vertices(\Psi)}{I_x = S}$ over $F(u)$. 
Hence 
\begin{equation*}
N_{\Gamma}(S) 
=
\charprod{M(\Gamma)^F}{\sgn^F}{H^F} 
 =
\charprod{M(\Psi)}{\sgn^F}{H^F}
= N_{\Psi}(S).
\end{equation*}

Now suppose $J \subseteq S$.  Let $\Psi_J$ be the $W_J$-graph
obtained from $\Psi$ by replacing
$I_x$ by $I_x \cap J$ for $x \in \vertices(\Psi)$.  
Also, let $\Gamma_J$ be the $W_J$-digraph obtained from 
$\Gamma$ by removing all edges with labels in
$S \setminus J$. 
Then
$M(\Gamma_J)^F \cong {M(\Gamma)^F \vert_{H_J^F}}
\cong {M(\Psi) \vert_{H_J^F}} \cong M(\Psi_J)$ as $H_J^F$-modules, 
and so the reasoning above gives
$N_{\Gamma_J}(J) = N_{\Psi_J}(J)$.
Therefore 
\[
\sum_{J \subseteq K \subseteq S} N_{\Gamma}(K)
=
N_{\Gamma_J}(J)
=
N_{\Psi_J}(J)
=
\sum_{J \subseteq K \subseteq S} N_{\Psi}(K).
\]
(Note $S$ itself is finite because $\edges(\Gamma)$ is
finite, so the sums above are finite.)
Thus 
part (i) of the theorem holds
by induction on $\abs{S \setminus J}$.

Let $\widehat{M}$ be the $H^F$-submodule of $M(\Psi)$ 
with basis
$\setof{x \in \vertices(\Psi)}{I_x \ne \emptyset}$.
By Lemma~\ref{lemma:indlemma},  
\[
M(\Psi)_{\ind^F} \cap \widehat{M} 
= 
\widehat{M}_{\ind^F}
=
\set{0}.
\]
Thus
\begin{equation*}
\begin{split}
\charprod{M(\Gamma)^F}{\ind^F}{H^F} 
 & =
\charprod{M(\Psi)}{\ind^F}{H^F} 
 =
\dim_{F(u)} M(\Psi)_{\ind^F} \\
& \le 
\dim_{F(u)} \left( M(\Psi) / \widehat{M} \right) 
 =
\abs{ \setof{ x \in \vertices(\Psi) } {I_x = \emptyset }} \\
& =
N_{\Psi}(\emptyset)
=
N_{\Gamma}(\emptyset),
\end{split}
\end{equation*}
with the last equality holding by part (i) of the theorem.  
Now, 
$\charprod{M(\Gamma)^F}{\ind^F}{H^F}$ is equal to the number of
connected components of $\Gamma$ by
Lemma~\ref{lemma:linearcharmults}(i), while
 $N_\Gamma(\emptyset)$ is equal to 
the number of sources of $\Gamma$. 
Therefore
$\Gamma$ has at least as many
sources as connected components. 
Because each connected component contains
at most one source by \cite{digraphpaper}, Theorem~1.5(i), it follows that 
every connected component of $\Gamma$  contains a 
(unique) source.  Hence every connected component of $\Gamma$
 is acyclic by \cite{digraphpaper}, Theorem~1.5(ii).
Therefore $\Gamma$ itself is acyclic, so part (ii) of the theorem holds and the 
proof of the theorem is complete.
\end{proof}


\section{An example}

In Figure~\ref{fig:affinecube}, a $W$-digraph $\Gamma$
is given for the affine group $W = W(\widetilde{A_2})$,
with generators $S=\set{r,s,t}$ satisfying
$n(r,s)=n(s,t)=n(r,t)=3$.  (The digraph $\Gamma$ is in 
fact a $W$-digraph by the classification given in Theorem~1.3
of \cite{digraphpaper}.)  Let $F$ be a subfield of $\complexes$.
We show $M(\Gamma)^F$ does not afford a $W$-graph over $F$,
arguing by contradiction.
\begin{figure}[H]
\caption{$W$-digraph $\Gamma$ for $W(\widetilde{A_2})$}
\centering
\input ./a2affine_cube.tex
\label{fig:affinecube}
\end{figure}

Suppose to the contrary that 
$\Psi=(\vertices(\Psi), I, \mu)$ is a $W$-graph
over $F \le \complexes$ such that $M(\Gamma)^F \cong M(\Psi)$.
Note that $\Gamma$ satisfies
$N_\Gamma(J)=1$ for all $J \subseteq S$.
Thus $N_\Psi(J)=1$ for $J \subseteq S$
by Theorem~\ref{theorem:wgraphacyclic}(i).
We can order $\vertices(\Psi)=\set{x_0,x_1,x_2,\dots,x_7}$ so
that, with $I_j = I_{x_j}$, we have
$$
\begin{aligned}
& \qquad I_0=  \emptyset,  I_1=\set{r}, I_2=\set{s}, I_3=\set{t}, \cr
& I_4=\set{r,s}, I_5=\set{r,t}, I_6=\set{s,t}, I_7=\set{r,s,t}. \cr
\end{aligned}
$$
Define  
$M_1 = \bigoplus_{1\le i \le 7} F(u) x_i$,
$M_2 = F(u) x_7$. It is clear that 
$M(\Psi) \ge M_1 \ge M_2$
as $H^F$-modules.  Moreover, 
$M(\Psi)/M_1 \cong M(\Psi)_{\ind^F}$ and
$M_2 = M(\Psi)_{\sgn^F}$ are 1-dimensional.
Put $\overline{M_1}= M_1 / M_2$, an $H^F$-module
with basis 
$$
\vertices_1=\set{x_1+M_2,x_2+M_2,\dots,x_6+M_2}
$$
over $F(u)$.
Put $\overline{x_j} = x_j + M_2$, and
define $I_1 :\vertices_1 \rightarrow \powerset(S)$
by $(I_1)_{\overline{x_j}} = I_{x_j}$, for $1\le j \le 6$.
Finally, define $\mu_1 : \vertices_1\times\vertices_1 \setminus \Delta
\rightarrow F$, by
$$
\mu_1(\overline{x_i},\overline{x_j})
=
\mu(x_i, x_j)
\qquad
\text{for $1 \le i, j \le 6, i \ne j$.}
$$
Then $\Psi_1=(\vertices_1,I_1,\mu_1)$ is a $W$-graph over $F$ 
with module 
$M(\Psi_1) = \overline{M_1}$.

Put $m_{ij}=\mu(x_i,x_j)$ for $1\le i,j \le 6$, $i \ne j$.
The matrices $B_r$, $B_s$, $B_t$ of $T_r$, $T_s$, $T_t$ 
acting on $M(\Psi_1)$ with respect to  the basis
$\vertices_1$ are 
$$
\begin{pmatrix}
-1 & m_{12}u & m_{13}u & 0 & 0 & m_{16}u \cr
0  &  u^2  &  0  &  0  &  0  &  0  \cr
0  &  0  &  u^2  &  0  &  0  &  0  \cr
0  &  m_{42}u & m_{43}u & -1 &  0 & m_{46}u \cr
0  &  m_{52}u & m_{53}u & 0  & -1  & m_{56}u \cr
0  &  0  &  0  &  0  &  0  &  u^2 \cr
\end{pmatrix},
\quad
\begin{pmatrix}
u^2  &  0  &  0  &  0  &  0  & 0 \cr
m_{21}u & -1 & m_{23}u & 0 & m_{25}u & 0 \cr
0 & 0 & u^2 & 0 & 0 & 0 \cr
m_{41}u & 0 & m_{43}u & -1 & m_{45} & 0 \cr
0  &  0  &  0  &  0  &  u^2  &  0 \cr
m_{61}u & 0 & m_{63}u & 0 & m_{65}u & -1 \cr
\end{pmatrix},
$$
and
$$ 
\begin{pmatrix}
u^2  &  0  &  0  &  0  &  0  & 0 \cr
0  &  u^2  &  0  &  0  &  0  &  0  \cr
m_{31}u & m_{32}u & -1 & m_{34}u & 0 & 0 \cr
0  &  0  &  0  &  u^2  &  0  &  0 \cr
m_{51}u & m_{52}u & 0 & m_{54}u & -1 & 0 \cr
m_{61}u & m_{62}u & 0 & m_{64}u & 0 & -1 \cr
\end{pmatrix},
$$
respectively.
For any finite sequence
$(s_1,s_2,\dots,s_\ell)$ in $S=\set{r,s,t}$, define
$$
T_{(s_1,s_2,\dots,s_\ell)}
=
T_{s_1}T_{s_2} \cdots T_{s_\ell}
\quad\text{and}\quad
B_{(s_1,s_2,\dots,s_\ell)}
=
B_{s_1}B_{s_2} \cdots B_{s_\ell}.
$$
Let $\chi_\Gamma : H^F \rightarrow F(u)$ be the character 
afforded by $M(\Gamma)^F$, and let
$\chi_\Psi$ and $\chi_{\Psi,1}$ denote the characters of
$M(\Psi)$ and $M(\Psi_1)$, respectively.  Then
$$
\chi_\Gamma = \chi_\Psi
= \chi_{\Psi,1} + \ind^F + \sgn^F.
$$
Thus
$$
\chi_\Gamma(T_{(s_1,s_2,\dots,s_\ell)}) 
- Tr(B_{(s_1,s_2,\dots,s_\ell)}) - u^{2\ell}-(-1)^\ell = 0
$$
for any finite sequence $(s_1,s_2,\dots,s_\ell)$ in $S$.

Direct calculations show
$$
\begin{aligned}
\chi_\Gamma(T_{(t, r, s, t, s, r)})  
- & Tr(B_{(t, r, s, t, s, r)}) - u^{12}-1 \cr
& = 
(1 - m_{56}m_{65}) u^2 + \cdots + (1 - m_{12} m_{21}) u^{10}, \cr
\end{aligned}
$$
where the omitted terms have degrees in $u$ between 3 and 9, inclusively.
This proves the first pair of equations in (\ref{eqn:rels1}).
\begin{equation}
m_{12}m_{21}=m_{56}m_{65}=1, \quad
m_{13}m_{31}=m_{46}m_{64}=1, \quad
m_{23}m_{32}=m_{45}m_{54}=1
\label{eqn:rels1}
\end{equation}
Similar calulations, using the sequences 
$(s, t, r, s, r, t)$ and $(r, s, t, r, t, s)$, establish the
remaining equations of (\ref{eqn:rels1}).
Calculations also show
$$
\begin{aligned}
\chi_\Gamma(T_{(s, r, s, r, t, s, r, t)})  
- & Tr(B_{(s, r, s, r, t, s, r, t)}) - u^{16}-1 \cr
& = 
(1 - m_{46}m_{54}m_{65} ) u^3 + \cdots + (1 - m_{12} m_{23}m_{31}) u^{13}, \cr
\end{aligned}
$$
so
\begin{equation}
m_{12} m_{23}m_{31} = m_{46}m_{54}m_{65} = 1.
\label{eqn:rels2}
\end{equation}
Next, the coefficient of $u^2$ in
$$
\chi_\Gamma(T_{(r, s, t, r, s)})  
- Tr(B_{(r, s, t, r, s)}) - u^{10}-1
$$
is 
$$
-2 + m_{16} m_{61} + m_{46} m_{64} + m_{56} m_{65},
$$
which is equal to $m_{16} m_{61}$ by (\ref{eqn:rels1}),
so 
the first equation of (\ref{eqn:rels3}) holds.
\begin{equation}
m_{16}m_{61}=0, \quad
m_{25}m_{52}=0, \quad
m_{34}m_{43}=0
\label{eqn:rels3}
\end{equation}
The remaining equations of (\ref{eqn:rels3}) are
verified by similar calculations using 
the sequences $(t, r, s, t, r)$ and
$(r, s, t, r, s)$.

Now, 
$$
\chi_\Gamma(T_{(r,s,t)})  - Tr(B_{(r,s,t)}) - u^6 + 1
=
\alpha u^2 - \beta u^3 + \gamma u^4,
$$
where
$$
\begin{cases}
\alpha = &  -3 + m_{34} m_{43} + m_{25} m_{52} + m_{45} m_{54}
 + m_{16} m_{61} + m_{46} m_{64} + m_{56} m_{65}, \cr
  \beta = & \ 
 m_{12} m_{23} m_{31} + m_{23} m_{34} m_{42} + m_{12} m_{25} m_{51} 
+ m_{25} m_{42} m_{54} \cr
& + m_{16} m_{31} m_{63} + 
  + m_{34} m_{46} m_{63} + m_{16} m_{51} m_{65} + m_{46} m_{54} m_{65}, \cr
  \gamma = & \ 
  3 - m_{12} m_{21} - m_{13} m_{31} - m_{23} m_{32} - m_{34} m_{43} - m_{25} m_{52} - m_{16} m_{61}. \cr
 \end{cases}
$$
Applying the relations 
(\ref{eqn:rels2}), it follows that
$$
\begin{aligned}
\beta = 
2 &+
 m_{23} m_{34} m_{42} + m_{12} m_{25} m_{51} + m_{25} m_{42} m_{54} \cr
&+ m_{16} m_{31} m_{63} + m_{34} m_{46} m_{63} + m_{16} m_{51} m_{65}. \cr
\end{aligned}
$$


We consider eight cases of the form $(a,b,c)=(0,0,0)$, where
$a \in \set{m_{16}, m_{61}}$, $b \in \set{m_{25}, m_{52}}$, 
$c \in \set{m_{34}, m_{43}}$.  These cases are exhaustive
by (\ref{eqn:rels3}). 
In each case it is shown that 
$\beta = 2 \ne 0$, giving a contradiction.


\vskip10pt
\noindent
Case~1.  $(m_{16}, m_{25}, m_{34})=(0,0,0)$.  In this case it is clear that
 $\beta=2$.


\vskip10pt
\noindent
Cases~2--7. $(a,b,c)=(0,0,0)$, 
where $(a,b,c)\ne (m_{16}, m_{25}, m_{34})$, 
$(a,b,c) \ne (m_{61}, m_{52}, m_{43})$.
In these six cases, we use the dihedral relations 
\begin{equation}
\begin{cases}
B_{(r,s,r)} - B_{(s,r,s)} &= {\bf 0} \cr
B_{(r,t,r)} - B_{(t,r,t)} &= {\bf 0} \cr
B_{(s,t,s)} - B_{(t,s,t)} &= {\bf 0} \cr
\end{cases}
\label{eqn:rels4}
\end{equation}
to show 
\begin{equation}
{\mathscr M} 
 = \set{0},
\label{eqn:rels5}
\end{equation}
where ${\mathscr M} = \set{m_{16},m_{61},m_{25},m_{52},m_{34},m_{43}}$.
In each case (\ref{eqn:rels5}) can be established by 
considering certain monomial entries in the matrices
on the left sides of (\ref{eqn:rels4}).

For example, consider the case $(m_{16},m_{25},m_{43})=(0,0,0)$.
In this case the $(6,1)$-entry of $B_{(r,t,r)} - B_{(t,r,t)}$ is
$$
-m_{61} (-1 + m_{13} m_{31} + m_{46} m_{64}) u^3,
$$
which is equal to
$-m_{61}u^3$ in view of (\ref{eqn:rels1}).
Therefore $m_{61}=0$.  Also, again using (\ref{eqn:rels1}),
the $(3,4)$- and  $(5,2)$-entries
of $B_{(s,t,s)} - B_{(t,s,t)}$ are equal to $-m_{34}u^3$ and
$-m_{52}u^3$, respectively, 
so $m_{34}=m_{52}=0$ and
(\ref{eqn:rels5}) holds.

As a second example, consider the case
$(m_{61},m_{25},m_{43})=(0,0,0)$.  Using (\ref{eqn:rels1}),
 the
$(1,6)$- and $(5,2)$-entries of $B_{(r,s,r)} - B_{(s,r,s)}$
are equal to $m_{16}u^3$ and $m_{52}u^3$, respectively, so
$m_{16}=m_{52}=0$.  Also, the $(3,4)$-entry
of $B_{(s,t,s)} - B_{(t,s,t)}$ is equal to $-m_{34}u^3$, and thus 
$m_{34}=0$, so again (\ref{eqn:rels5}) holds.

Similary calculations establish (\ref{eqn:rels5}) in the remaining
cases in this group.  
From (\ref{eqn:rels5}), it follows that $\beta=2$.  
(In any of Cases 2--7, to show $m_{ij}=0$, where 
$m_{ij} \in {\mathscr M} \setminus \set{a,b,c}$, it suffices to 
look at the $(i,j)$-entry of
one matrix on the left side of (\ref{eqn:rels4}) and apply 
(\ref{eqn:rels1}).  The author has no explanation for this 
pattern.)

\vskip10pt
\noindent
Case~8. $(m_{61},m_{52},m_{43})=(0,0,0)$. 
The $(1,5)$-entry of $B_{(r,s,r)} - B_{(s,r,s)}$ is
$$
-(m_{12} m_{25} + m_{16} m_{65}) (u^2 + u^4), 
$$
and thus 
$m_{12} m_{25} + m_{16} m_{65} = 0$.
Also, the $(3,6)$-entry of $B_{(r,t,r)} - B_{(t,r,t)}$ is
$$
(m_{16} m_{31} + m_{34} m_{46}) (u^2 + u^4),
$$ 
so
$m_{16} m_{31} + m_{34} m_{46} = 0$.
Finally, the $(2,4)$-entry of $B_{(s,t,s)} - B_{(t,s,t)}$ is
$$
-(m_{23}  m_{34}  + m_{25}  m_{54} ) (u^2 + u^4),
$$
and thus
$m_{23}  m_{34}  + m_{25}  m_{54} = 0$. 
Hence
$$
\begin{aligned}
\beta & = 
2+ (m_{23}  m_{34}  + m_{25}  m_{54}) m_{42} +
   (m_{16} m_{31} 
   + m_{34} m_{46}) m_{63} \cr
&   \qquad + (m_{12} m_{25} + m_{16} m_{65}) m_{51} = 2. \cr
 \end{aligned}
 $$
 
 Since $\beta=2$ in all cases, 
 we have arrived at a contradiction Therefore $M(\Gamma)^F$ is
 not isomorphic to $M(\Psi)$.


\bibliographystyle{plain}
\bibliography{combined}


\vfil\eject
\enddocument
\bye

%% file: header.tex
\usepackage[paperwidth=8.5in, textheight=9.5in]{geometry}
\usepackage{amsmath,amsthm,amssymb}
\usepackage{tikz}
\usetikzlibrary{arrows,decorations.pathmorphing,backgrounds,positioning,fit,petri}
\usepackage{enumerate,array,float,verbatim}
\usepackage[labelsep=space]{caption}
\usepackage{mathrsfs}

\theoremstyle{plain}
\newtheorem{theorem}{Theorem}[section]
\newtheorem{lemma}[theorem]{Lemma}

\theoremstyle{definition}

\numberwithin{subcase}{case}
\numberwithin{subsubcase}{subcase}
\numberwithin{table}{section}
\numberwithin{equation}{section}

\newcommand{\singleedge}[5]{\begin{tikzpicture}%
[node distance=0.8cm,baseline={(left.base)},%
pre/.style={<-,shorten <=1pt,>=angle 45},%
post/.style={->,shorten >=1pt,>=angle 45}],%
rectangle/.style={inner sep=0pt,minimum size=4mm}];%
\node[rectangle] (left) {$\strut$#1};%
\node[rectangle] (right) [right=of left] {$\strut$#2}%
edge [#4,#5] node[yshift=5pt] {#3} (left);%
\end{tikzpicture}}
\newcommand{\solidedge}[3]{\singleedge{#1}{#2}{#3}{pre}{solid}}
\newcommand{\dashededge}[3]{\singleedge{#1}{#2}{#3}{pre}{dashed}}

\newcommand{\infinity}{\infty}
\newcommand{\set}[1]{\left\{#1\right\}}
\newcommand{\setof}[2]{\left\{{#1}\mid{#2}\right\}}

\newcommand{\abs}[1]{\left\vert{#1}\right\vert}
\newcommand{\rationals}{\mathbb Q}%
\newcommand{\complexes}{\mathbb C}

\newcommand{\charprod}[3]{\left<{#1},\strut{#2}\right>_{#3}}

\DeclareMathOperator{\ind}{ind}
\DeclareMathOperator{\sgn}{sgn}

\DeclareMathOperator{\vertices}{{\mathscr V}}
\DeclareMathOperator{\edges}{{\mathscr E}}
\DeclareMathOperator{\powerset}{{\mathscr P}}
\DeclareMathOperator{\In}{In}


%% file: a2affine_cube.tex
\begin{tikzpicture}
[node distance=1.5cm,%
pre/.style={<-,shorten <=1pt,>=angle 45},%
post/.style={->,shorten >=1pt,>=angle 45}];%
rectangle/.style={inner sep=0pt,minimum size=5mm}];%
\tikzstyle{mydot}=[circle,on grid,draw=black,fill=black,inner sep=0pt,minimum size=1mm];
\node[rectangle] (gamma5)								{${\strut}\gamma_5$};
\node[rectangle] (gamma6) [right=of gamma5,xshift=30mm]		{${\strut}\gamma_6$}
	edge[pre,dashed]	node[yshift=2.5mm]			{$r$}		(gamma5);
\node[rectangle] (gamma1) [below=of gamma5,xshift=15mm,yshift=10mm] {${\strut}\gamma_1$}
	edge[post]			node[xshift=-1mm,yshift=-2mm]	{$t$}		(gamma5);
\node[rectangle] (gamma2) [below=of gamma6,xshift=-15mm,yshift=10mm] {${\strut}\gamma_2$}
	edge[pre]				node[yshift=2.5mm]		{$r$}		(gamma1)
	edge[post,dashed]	node[xshift=1mm,yshift=-2mm]	{$t$}		(gamma6);
\node[rectangle] (gamma3) [below=of gamma1,yshift=1mm]		{${\strut}\gamma_3$}
	edge[pre]				node[xshift=-2.5mm]		{$s$}		(gamma1);
\node[rectangle] (gamma4) [below=of gamma2,yshift=1mm]		{${\strut}\gamma_4$}
	edge[pre,dashed]		node[yshift=2.5mm]		{$r$}		(gamma3)
	edge[pre,dashed]		node[xshift=2.5mm]		{$s$}		(gamma2);
\node[rectangle] (gamma7) [below=of gamma3,xshift=-15mm,yshift=10mm] {${\strut}\gamma_7$}
	edge[pre,dashed]		node[xshift=-2.5mm]		{$s$}		(gamma5)
	edge[pre,dashed]		node[xshift=-1mm,yshift=2mm] {$t$} (gamma3);
\node[rectangle] (gamma8) [below=of gamma4,xshift=15mm,yshift=10mm] {${\strut}\gamma_8$}
	edge[pre]				node[yshift=2.5mm]		{$r$}		(gamma7)
	edge[pre]			node[xshift=1mm,yshift=2mm] {$t$}		(gamma4)
	edge[pre]				node[xshift=2.5mm]		{$s$}		(gamma6);
\end{tikzpicture}

%% file: wdgrii.bbl
\begin{thebibliography}{1}

\bibitem{digraphpaper}
D.~Alvis.
\newblock A {C}lass of {R}epresentations of {H}ecke {A}lgebras.
\newblock {\em Bull. Inst. Math. Acad. Sinica (N.S.)}, 11(2):301--342, 2016.

\bibitem{gyojawgraph}
Akihiko Gyoja.
\newblock On the existence of a ${W}$-graph for an irreducible representation
  of a {C}oxeter group.
\newblock {\em J.~Alg.}, 86:422--438, 1984.

\end{thebibliography}
